\documentclass[a4paper,11pt]{article}
\pdfoutput=1

\usepackage[a4paper,tmargin=3truecm,bmargin=3truecm,rmargin=2.5truecm,
lmargin=2.5truecm,twoside,verbose=true]{geometry}
\usepackage{cancel,graphicx}

\usepackage{amsmath,amssymb}
\usepackage[amsmath, hyperref, thmmarks]{ntheorem}
\usepackage[all]{xypic}
\usepackage[pdftex]{hyperref}
\usepackage[english]{babel}

\numberwithin{equation}{section}

\allowdisplaybreaks[1]

\def\question#1{}

\newcommand\caB{{\mathcal B}}
\newcommand\caC{{\mathcal C}}
\newcommand\caD{{\mathcal D}}

\newcommand\caL{{\mathcal L}}

\newcommand\wx{{\widetilde x}}
\newcommand\gone{{ \mathchoice {1\mskip-4mu\mathrm{l} } {1\mskip-4mu\mathrm{l} }{1\mskip-4.5mu\mathrm{l} } {1\mskip-5mu\mathrm{l}} }}
\newcommand\gR{{\mathbb R}}

\newcommand\gC{{\mathbb C}}

\newcommand\gB{{\mathbb B}}
\newcommand\gN{{\mathbb N}}
\newcommand\gZ{{\mathbb Z}}

\newcommand\algzero{{\mathsf 0}}
\newcommand\superA{{\mathcal A}}
\newcommand\algA{{\mathbf A}}

\newcommand\algB{{\mathbf B}}

\newcommand\ehH{\mathcal H}

\newcommand\eps{{\varepsilon}}

\newcommand\fois{\mathord{\cdot}}
\DeclareMathOperator{\Tr}{Tr}

\newcommand\dd{{\text{\textup{d}}}}

\newcommand\norm{\mathord{\parallel}}

\theoremsymbol{}
\theorembodyfont{\slshape}
\theoremheaderfont{\normalfont\bfseries}
\theoremseparator{}
\newtheorem{Theorem}{Theorem}[section]
\newtheorem{theorem}[Theorem]{Theorem}

\newtheorem{proposition}[Theorem]{Proposition}

\theorembodyfont{\upshape}
\theoremsymbol{\ensuremath{\blacklozenge}}

\newtheorem{definition}[Theorem]{Definition}

\theoremstyle{nonumberplain}
\theoremheaderfont{\scshape}
\theorembodyfont{\normalfont}
\theoremsymbol{\ensuremath{\blacksquare}}

\qedsymbol{\ensuremath{_\blacksquare}}
\theoremclass{LaTeX}

\title{Non-formal deformation quantization of abelian supergroups\footnote{Work
supported by the Belgian Interuniversity Attraction Pole (IAP) within the framework ``Nonlinear systems, stochastic processes, and statistical mechanics'' (NOSY).}}
\date{}
\author{Axel de Goursac}

\begin{document}

\maketitle
\vspace*{-1cm}
\begin{center}
\textit{D\'epartement de Math\'ematiques, Universit\'e Catholique de Louvain,\\ Chemin du Cyclotron 2, 1348 Louvain-la-Neuve, Belgium\\
e-mail: \texttt{axelmg@melix.net}}\\
\end{center}%

\vskip 2cm

\begin{abstract}
We review recent works concerning deformation quantization of abelian supergroups. Indeed, we expose the construction of an induced representation of the Heisenberg supergroup and an associated pseudodifferential calculus by using Kirillov's orbits method. Then, a star-product is built on the abelian supergroup $\gR^{m|n}$ together with a universal deformation formula for its actions. Using topological Hopf algebras, we reformulate this deformation as a continuous twist on comodule-algebras. We also introduce the notion of C*-superalgebra which is natural and compatible with the deformation. Finally, we show some applications in supergeometry and theoretical physics.
\end{abstract}

\vskip 0.5cm

\hspace*{1cm} induced representation; Hopf algebra; renormalization

\vfill

\pagebreak

\section{Introduction}

\subsection{Motivations}

This paper is a survey of recent works \cite{Bieliavsky:2010su,deGoursac:2011kv,deGoursac:2010zb} on the non-formal deformation quantization of abelian supergroups and its applications in theoretical physics. These works are set in the perspective of the construction of a noncommutative supergeometry, namely a noncommutative geometry \cite{Connes:1994} for graded spaces (see \cite{deGoursac:2008bd} for an algebraic approach of this construction). The objects studied here are deformation quantization \cite{Bayen:1978,Kontsevich:2003} of abelian supergroups, consisting in the data of a noncommutative product on the algebra of smooth functions $C^\infty(M)$ and depending on a deformation parameter. For more about formal deformations of groups see \cite{Drinfeld:1989,Giaquinto:1998}, and \cite{Rieffel:1989,Bieliavsky:2002,Bieliavsky:2007} in the non-formal case.

The paper is organized as follows. After having recalled the necessary supergeometrical framework (subsection \ref{subsec-setting}), we briefly review the different mathematical domains where the works \cite{Bieliavsky:2010su,deGoursac:2011kv} can be applied for the noncommutative supergeometry.
\begin{itemize}
\item In harmonic analysis with the construction of a unitary induced representation of the Heisenberg supergroup and a pseudodifferential calculus (subsections \ref{subsec-unindrepr} and \ref{subsec-qu}).
\item In deformation quantization with the construction of a non-formal star-product on the abelian supergroup $\gR^{m|n}$ (subsection \ref{subsec-starprod}) as well as a universal deformation formula (subsection \ref{subsec-udf}).
\item In topological Hopf algebras with the reformulation of the above deformation as a continuous twist on Fr\'echet comodule-algebras (subsection \ref{subsec-twist}).
\item In operator algebras with the definition of Hilbert superspaces and C*-superalgebras, which are appropriate notions for this setting of deformation of supergroups (subsections \ref{subsec-oper} and \ref{subsec-udfcstar}).
\item In supergeometry where the previous constructions can be used to deform other supermanifolds (subsection \ref{subsec-supergeom}).
\item In theoretical physics where a renormalizable noncommutative quantum field theory (QFT) is interpreted by using a deformation of a superspace (subsection \ref{subsec-qft}).
\end{itemize}

\subsection{Setting}
\label{subsec-setting}

We will work essentially in the concrete setting of supergeometry developed by DeWitt, Rogers, Tuynman,... \cite{DeWitt:1984,Rogers:2007,Tuynman:2005}. Let $\superA=\bigwedge V$, where $V$ is an real infinite-dimensional vector space, it is a $\gZ_2$-graded commutative algebra: $\superA=\superA_0\oplus\superA_1$ satisfying:
\begin{equation*}
\forall a,b\in\superA\quad:\quad ab=(-1)^{|a||b|}ba,
\end{equation*}
where $|a|\in\gZ_2$ denotes the degree of the homogeneous element $a$. Then, we define the superspace of dimension $m|n$ to be $\gR^{m|n}=(\superA_0)^m\times(\superA_1)^n$. A smooth function on such a superspace $f\in C^\infty(\gR^{m|n})$ can be decomposed as
\begin{equation}
f(x,\xi)=\sum_{I} f_I(x)\xi^I,\label{eq-decompsuper}
\end{equation}
for $x\in\gR^m$, $\xi\in\gR^{0|n}$ (anticommuting coordinates), $\xi^I:=\prod_{i\in I}\xi^i$, and where $I$ is summed over the (ordered) subsets of $\{1,\dots,n\}$ and where $f_I\in C^\infty(\gR^m)$. We also denote by $\gB:\gR^{m|n}\to\gR^m$ the body map. For more details, in particular about the construction of supermanifolds, see e.g. \cite{Rogers:2007,Tuynman:2005}.

In the following, we will consider complex-valued superfunctions. Let $\superA_\gC=\superA\otimes\gC$. One can then define a complex conjugation as
\begin{equation*}
\forall a\in\superA,\,\forall\lambda\in\gC\quad:\quad \overline{a\otimes\lambda}:=a\otimes\overline{\lambda},
\end{equation*}
so that the property $\overline{ab}=(-1)^{|a||b|}\overline{b}\overline{a}$, $\forall a,b\in\superA_\gC$ stands. By using the decomposition \eqref{eq-decompsuper}, the Lebesgue-Berezin integral is well-defined on $\gR^{m|n}$:
\begin{equation*}
\int_{\gR^{m|n}}\dd z\, f(z)=\int_{\gR^m}\dd x\,f_{\{1,\dots,n\}}(x),
\end{equation*}
as well as the Hodge operation \cite{Bieliavsky:2010su}:
\begin{equation}
(\ast f)(x,\xi)=\sum_I \eps(I,\complement I) f_I(x) \xi^{\complement I},\label{eq-hodge}
\end{equation}
where $\eps(I,J)$ is the signature of the permutation needed to order the set $I\cup J$, and $\complement I$ is the complementary subset of $I$ in $\{1,\dots,n\}$. Two scalar products have then been constructed on the complex superfunctions $C^\infty(\gR^{m|n})$:
\begin{align}
&\langle f,g\rangle=\int\dd x\,\dd\xi\, \overline{f(x,\xi)}g(x,\xi)=\sum_I\eps(I,\complement I)\int\dd x\overline{f_I(x)}g_{\complement I}(x),\nonumber\\
&\left(f,g\right)=\int\dd x\,\dd\xi\, \overline{f(x,\xi)}(\ast g)(x,\xi)=\sum_I\int\dd x\overline{f_I(x)}g_{I}(x).\label{eq-func-scalprod}
\end{align}
The first one is superhermitian whereas the second one is hermitian positive definite. The completion of $C^\infty(\gR^{m|n})$ under the second scalar product is denoted by $L^2(\gR^{m|n})$.
\medskip

We also recall the definition of the Heisenberg supergroup $G$. As a supermanifold, it is isomorphic to $\gR^{m|n}\times\gR^{1|0}$, where $m$ is an even number. Let $\omega$ be the even symplectic structure on $\gR^{m|n}\subset G$ associated to the square matrix $\begin{pmatrix} \omega_0 & 0 \\ 0 & 2\gone \end{pmatrix}$ of size $m+n$ in the homogeneous canonical basis, with $\omega_0=\begin{pmatrix} 0 & \gone \\ -\gone & 0 \end{pmatrix}$ of size $m$. The smooth group law of $G$ is then given by
\begin{equation*}
\forall (x,\xi),(y,\eta)\in\gR^{m|n},\,\forall a,b\in\gR^{1|0}\quad:\quad (x,\xi,a)\fois (y,\eta,b)=(x+y,\xi+\eta,a+b+\frac 12\omega_0(x,y)+\xi\eta).
\end{equation*}
$G$ is a non-abelian supergroup of dimension $m+1|n$, with neutral element $(0,0,0)$, and inverse $(x,\xi,a)^{-1}=(-x,-\xi,-a)$.

\section{Construction of the deformation quantization}

\subsection{Unitary induced representation}
\label{subsec-unindrepr}

To construct a unitary induced representation of the Heisenberg supergroup $G$, we use the Kirillov's orbits method \cite{Kirillov:1976}. We choose a real polarization associated to a coadjoint orbit of $G$. The corresponding normal supergroup can be written as $B=W\times K$, where $W\simeq \gR^{\frac m2|0}$ is a maximal isotropic subspace of $\gR^{m|n}\subset G$ and $K=\gR^{1|0}$. The associated character $\chi:B\to\superA_\gC$ can be expressed ad $\chi(w,a)=e^{\frac{i}{\theta}a}$ with $\theta$ a real non-zero parameter of the orbit. If $Q\simeq\gR^{\frac m2|n}$ is the dual of $W$ for the symplectic form $\omega$, we have the decomposition
\begin{equation}
G\simeq Q\times B\label{eq-decompG}.
\end{equation}
We consider now the left regular action $g\in G\mapsto L^\ast_{g^{-1}}\hat\varphi$ of $G$ on $B$-equivariant superfunctions $\hat\varphi\in C^\infty(G)^B$, that are smooth superfunctions $\hat\varphi$ satisfying
\begin{equation*}
\forall g\in G,\ \forall b\in B\quad:\quad \hat\varphi(g\fois b)=\chi(b^{-1})\hat\varphi(g).
\end{equation*}
Because of \eqref{eq-decompG}, we have the isomorphism $C^{\infty}(G)^B\simeq C^\infty(Q)$.
\begin{proposition}[Induced representation]
The reexpression of the left regular action on $C^\infty(Q)$ gives a representation $U:G\to\caL_\superA(\caD(Q)\otimes\superA)$ defined by $\forall\varphi\in C^\infty(Q)$, $\forall (x,\xi,w,a)\in G$, $\forall(x_0,\xi_0)\in Q$,
\begin{equation*}
U(x,\xi,w,a)\varphi(x_0,\xi_0)= e^{\frac{i}{\theta}(a+\frac12\omega_0(x-2x_0,w)+\xi\xi_0)} \varphi(x_0-x,\xi_0-\xi),
\end{equation*}
where $\caD(Q)$ is the space of smooth superfunctions with compact support on $Q$, and $\caL_\superA$ is the space of $\superA$-linear endomorphisms. $U$ is the representation of $G$ induced from $(B,\chi)$, which generalizes the Schr\"odinger representation to the graded setting.
\end{proposition}

Note that there are several possible choices for the scalar product on $C^\infty(Q)$. It turns out that the representation we have built is unitary for the superhermitian scalar product of Equation \eqref{eq-func-scalprod}
\begin{equation*}
\forall g\in G\quad :\quad \langle U(g)\varphi,U(g)\psi\rangle=\langle\varphi,\psi\rangle,
\end{equation*}
which seems to be natural in this framework (see \cite{Alldridge:2011iz} for another choice). This choice of scalar product is indeed determinant of the harmonic analysis we can construct on the Heisenberg supergroup.

\subsection{Quantization map}
\label{subsec-qu}

Using this unitary induced representation $U$, we will construct now a quantization map $\Omega$. Let us introduce the map $\Sigma\in \caL_\superA(\caD(Q)\otimes\superA)$ given by
\begin{equation*}
\forall (x_0,\xi_0)\in Q\quad:\quad (\Sigma\varphi)(x_0,\xi_0)=\gamma\int\dd\xi_1\ e^{-\frac{i\alpha}{\theta}\xi_1\xi_0}\varphi(-x_0,\xi_1),
\end{equation*}
where $\alpha$ is another real parameter of the deformation (different from $-1$ and $0$), and $\gamma=\frac{(-1)^{\frac{n(n+1)}{2}}(i\theta)^n}{(\pi\theta)^{\frac m2}(1+\alpha)^n}$. By setting $r=\frac{(-1)^n\alpha^n}{(\pi\theta)^{\frac m2}(1+\alpha)^n}\gamma$, we see that $\Sigma$ is almost-involutive and almost-unitary
\begin{equation}
\Sigma^2=r\gone,\quad \Sigma^\dag=r\Sigma.\label{eq-propsigma}
\end{equation}
The quantization map $\Omega:G\to \caL_\superA(\caD(Q)\otimes\superA)$ defined by
\begin{equation*}
\forall g\in G\quad:\quad \omega(g)=U(g)\Sigma U(g^{-1}),
\end{equation*}
is constant on the left classes of $G$ by $K$, and is well-defined on the quotient $G/K=Q\times W\simeq\gR^{m|n}$. Attention, $\gR^{m|n}$ is an abelian supergroup but $\Omega$ is not a group homomorphism.

We define the quantization map on superfunctions by
\begin{equation*}
\forall f\in L^1(\gR^{m|n})\quad:\quad \Omega(f)=\int_{\gR^{m|n}}\dd z\ f(z)\Omega(z).
\end{equation*}
Its explicit expression is given by
\begin{multline}
(\Omega(f)\varphi)(x_0,\xi_0)=\gamma\int_M\dd x\dd\xi\dd w f(x,\xi,w)\\
\int\dd\xi_1\ e^{\frac{i}{\theta}(2\omega_0(x-x_0,w)+\xi\xi_0-\alpha\xi_1\xi_0 -(\alpha+1)\xi\xi_1)}
\varphi(2x-x_0,\xi+\xi_1).\label{eq-omf}
\end{multline}
Moreover, this map can be continuously extended to $\Omega:L^1(\gR^{m|n})\to\caL(L^2(Q))$, where $\caL$ means the bounded operators or continuous linear maps. It is an extension of the Weyl quantization to the graded setting.

We want to enlarge the quantization map to the space $\caB(\gR^{m|n})$ of smooth superfunctions whose every derivative are bounded. The expression \eqref{eq-omf} is meaningless for non-integrable functions $f$, so that we consider the notion of oscillating integral as follows. Let us introduce the operator $O$:
\begin{equation*}
(O\fois f)(x,\xi,w)=(1-\Delta_{(x,w)}) \Big(\frac{1}{1+x^2+w^2}f(x,\xi,w)\Big)\ ,
\end{equation*}
for a smooth superfunction $f$ with compact support and where $\Delta_{(x,w)}$ denotes the Laplacian with respect to the variables $(x,w)$. With an integration by parts, we obtain
\begin{equation}
\int\dd x\dd\xi\dd w\ e^{i\omega_0(x,w)}f(x,\xi,w)= \int\dd x\dd\xi\dd w\ e^{i\omega_0(x,w)}(O^k\fois f)(x,\xi,w),\label{eq-oscil}
\end{equation}
for any $k\in\gN$. Moreover, one can also prove that there exist functions $b^\alpha\in\caB(\gR^m)$ such that
\begin{equation}
(O^k\fois f)(x,\xi,w)=\frac{1}{(1+x^2+w^2)^k}\sum_{\alpha\in\gN^m,\,|\alpha|\leq 2k} b^\alpha(x,w) D^\alpha f(x,\xi,w).\label{eq-oscil2}
\end{equation}
Consequently, there exists a number $k$ such that $(O^k\fois f)\in L^1(\gR^{m|n})$ for any $f\in\caB(\gR^{m|n})$. In the following, we mean by oscillating integral of $f$ the RHS member of \eqref{eq-oscil}.
\begin{proposition}[Quantization map]
\label{prop-qmap}
The quantization map can now be extended into the continuous map
\begin{equation*}
\Omega:\caB(\gR^{m|n})\to\caL(L^2(Q)),
\end{equation*}
where the expression \eqref{eq-omf} is understood as an oscillating integral.
\end{proposition}
If $f$ is in a symbol space, Equation \eqref{eq-omf} provides that $\Omega(f)$ is a pseudodifferential operator on $Q$. In a more general way, a symbol calculus $\caB^\mu(\gR^{m|n})$ has been introduced in \cite{Bieliavsky:2010su} which defines via $\Omega$ a pseudodifferential calculus.

\subsection{Star-product}
\label{subsec-starprod}

In this subsection, we express the (non-formal) star-product on $\gR^{m|n}$ corresponding to the composition of operators via this quantization map $\Omega$. It is defined by the following formula: for $(x,\xi,w)\in\gR^{m|n}$,
\begin{multline}
(f_1\star f_2)(x,\xi,w)=\kappa\int\dd x_1\dd\xi_1\dd w_1\dd x_2\dd\xi_2\dd w_2\ f_1(x_1,\xi_1,w_1)f_2(x_2,\xi_2,w_2)\\
e^{\frac{2i}{\theta}(\omega_0(x_2-x,w_1)+\omega_0(x-x_1,w_2)+\omega_0(x_1-x_2,w))}
e^{-i\frac{(1+\alpha)^2}{\alpha\theta}(\xi\xi_1+\xi_1\xi_2+\xi_2\xi)},\label{eq-moy}
\end{multline}
where $\kappa=\frac{(i\theta)^n(-1)^{\frac{n(n+1)}{2}}\alpha^n}{(\pi\theta)^m(1+\alpha)^{2n}}$, and $\theta$, $\alpha$ are the parameters of the deformation.

\begin{theorem}[Star-product]
It turns the space $\caB(\gR^{m|n})$ into an associative (but noncommutative) algebra and satisfies:
\begin{equation*}
\forall f_1,f_2\in\caB(\gR^{m|n})\quad:\quad \Omega(f_1\star f_2)=\Omega(f_1)\Omega(f_2),\qquad \overline{f_1\star f_2}=(-1)^{|f_1||f_2|}\overline{f_2}\star\overline{f_1}.
\end{equation*}
It can also be shown that the tracial property holds. For $f_1,f_2\in L^2(\gR^{m|n})$, $f_1\star f_2$ is well-defined, belongs to $L^1(\gR^{m|n})$ and satisfies:
\begin{equation*}
\int\dd z\, (f_1\star f_2)(z)=\int \dd z\, f_1(z)f_2(z).
\end{equation*}
\end{theorem}
The algebra $\caB(\gR^{m|n})$ endowed with this deformed product is then isomorphic to the tensor product of the Moyal algebra $\gR^m_\theta$ with the Clifford algebra $Cl(n,\gC)$.

Let us consider a rescaling of the symplectic form $\omega$: $\tilde\omega=\begin{pmatrix} \omega_0 & 0 \\ 0 & \frac{(1+\alpha)^2}{2\alpha}\gone \end{pmatrix}$. Then, the commutative limit can be computed: for $f_1,f_2\in \caB(\gR^{m|n})$,
\begin{equation*}
f_1\star f_2|_{\theta=0}=f_1 f_2,\qquad \frac{1}{\theta}(f_1\star f_2-(-1)^{|f_1||f_2|}f_2\star f_1)|_{\theta=0} =\frac{-i(1+\alpha)^{n^2-2n}}{\alpha^n}(-1)^{|f_1||\mu|} \tilde\omega^{-1}_{\nu\mu}\partial_\mu f_1\partial_\nu f_2,
\end{equation*}
where $\mu,\nu\in\{1,\dots,m+n\}$ and $|\mu|\in\gZ_2$ denotes the degree of the coordinate $z_\mu$. This means that the star-product \eqref{eq-moy} is a deformation quantization of the symplectic supermanifold $(\gR^{m|n},\tilde\omega^{-1})$ in the direction of its standard Poisson bracket.

The internal symmetries of a star-product are the diffeomorphisms $\phi$ which let the star-product invariant: $\phi^\ast(f_1\star f_2)=(\phi^\ast f_1)\star(\phi^\ast f_2)$.
\begin{proposition}[Symmetries]
The internal symmetry group of the star-product \eqref{eq-moy} is the semidirect product of the orthosymplectic group of $(\gR^{m|n},\tilde\omega^{-1})$ with the translations: $OSP(n,\frac m2)\ltimes \gR^{m|n}$.
\end{proposition}

\section{Universal deformation formula}

\subsection{UDF for Fr\'echet algebras}
\label{subsec-udf}

In this subsection, we study how the star-product we have constructed above can also deform the associative structure of any Fr\'echet algebra on which the abelian supergroup $\gR^{m|n}$ acts. Let $\rho$ be an action of $\gR^{m|n}$ on a Fr\'echet algebra $(\algA,|\fois|_j)$:
\begin{equation*}
\rho:\gR^{m|n}\times(\algA\otimes\superA)\to(\algA\otimes\superA),
\end{equation*}
We assume that $\rho$ satisfies the axioms
\begin{itemize}
\item $\rho_0=\text{id}$; $\forall z_1,z_2\in\gR^{m|n}$, $\rho_{z_1+z_2}=\rho_{z_1}\rho_{z_2}$.
\item $\forall z\in\gR^{m|n}$, $\rho_z:(\algA\otimes\superA)\to(\algA\otimes\superA)$ is an $\superA$-linear automorphism of algebras.
\item If we write $z=(y,\xi)\in\gR^{m|n}$, then the action can be expanded as: $\rho_{(y,\xi)}(a)=\sum_I \rho_y(a)_I\xi^I$, and $\forall a\in\algA$, $\forall I$, $y\mapsto \rho_y(a)_I$ is $\algA$-valued and continuous.
\item There exists a constant $C>0$ such that
\begin{equation*}
\forall a\in\algA,\, \forall I,\,\forall j,\quad \exists k,\quad \forall y\in\gB M, \quad |\rho_y(a)_I|_j\leq C|a|_k.
\end{equation*}
\end{itemize}
The expression \eqref{eq-moy} defines also (by taking the product of $f_1$ and $f_2$ in $\algA$) a star-product on the Fr\'echet space $\caB_\algA(\gR^{m|n})$ of $\algA$-valued smooth superfunctions with every derivative bounded, endowed with the seminorms $|f|_{j,\alpha} =\sup_{x\in\gR^m}\{\sum_I|D^\alpha f_I(x)|_j\}$.

If we denote by $\algA^\infty$ the smooth vectors \cite{Warner:1972}, i.e. the dense subset of elements $a$ of $\algA$ such that the map $\rho^a:z\mapsto\rho_z(a)$ is smooth on $\gR^{m|n}$, then it can be shown that
\begin{equation}
\forall a\in\algA^\infty\quad:\quad \rho^a\in\caB_{\algA}(\gR^{m|n})\label{eq-rhobound}
\end{equation}

\begin{theorem}[Universal Deformation Formula]
\label{thm-udf}
Consequently, the expression $a\star_\rho b:=(\rho^a\star\rho^b)(0)$ for smooth vectors $a,b\in\algA^\infty$ is well-defined. Endowed with the seminorms $|a|_{j,\alpha}:=|\rho^a|_{j,\alpha}$, 
$(\algA^\infty,\star_\rho)$ is a Fr\'echet algebra.
\end{theorem}

We have indeed constructed a deformation of the associative $\algA$, or more precisely of a dense subalgebra $\algA^\infty$, by using the star-product \eqref{eq-moy} and the action $\rho$ of $\gR^{m|n}$ on $\algA$. The examples of such deformations are numerous. We will see some of them closely related to supergeometry in subsection \ref{subsec-supergeom}.

\subsection{Twist on comodule-algebras}
\label{subsec-twist}

Let us denote the Fr\'echet algebra $H:=(\caB(\gR^{m|n}),|\fois|_\alpha)$. We show in this subsection (see \cite{deGoursac:2011kv} for the proofs) that $H$ has a structure of Fr\'echet-Hopf algebra and that the universal deformation formula can be reformulated into a continuous twist on the comodule-algebras of $H$.

First, we recall the definition of two topological tensor products. Let $(\algA,|\fois|_j)$ and $(\algB,|\fois|_k)$ be two arbitrary Fr\'echet algebras. The projective tensor product \cite{Grothendieck:1966} $\algA\widehat\otimes_\pi\algB$ is the completion of the algebraic tensor product $\algA\otimes\algB$ for the family of seminorms: $\forall c\in\algA\otimes\algB$,
\begin{equation*}
\pi_{j,k}(c)=\inf\{|\sum_i |a_i|_j|b_i|_k,\quad c=\sum_i a_i\otimes b_i\},
\end{equation*}
where the infimum is taken over all decompositions $c=\sum_i a_i\otimes b_i$. The topological tensor product of $(\algA,|\fois|_j)$ with $H$, denoted by $\algA\widehat\otimes_{\tau}H$, is defined as the completion of the algebraic tensor product for the family of seminorms of $\caB_\algA(\gR^{m|n})$:
\begin{equation*}
\tau_{j,\alpha}(f)=|f|_{j,\alpha} =\sup_{x\in\gR^m}\{\sum_I|D^\alpha f_I(x)|_j\}.
\end{equation*}
One can show that $H\widehat\otimes_{\tau}H\simeq \caB(\gR^{m|n}\times \gR^{m|n})$, and by definition, $\algA\widehat\otimes_{\tau}H\simeq \caB_\algA(\gR^{m|n})$. Since $H$ is not nuclear, $\tau$ is not equivalent to $\pi$ in general. However, they are compatible in the following sense: the exchange map
\begin{equation*}
\sigma_{23}:(\algA\widehat\otimes_\tau H)\widehat\otimes_\pi(\algA\widehat\otimes_\tau H)\to (\algA\widehat\otimes_\pi\algA)\widehat\otimes_\tau(H\widehat\otimes_\tau H),
\end{equation*}
defined by $\sigma_{23}(a_1\otimes f_1\otimes a_2\otimes f_2)=a_1\otimes a_2\otimes f_1\otimes f_2$, is continuous, for any Fr\'echet algebra $\algA$. This is essential for Proposition \ref{prop-comod}.

The structure of Hopf algebra we expose here is natural for an algebra of functions on a supergroup.
\begin{itemize}
\item the product $\tilde\mu:H\widehat\otimes_\tau H\to H$ by $\tilde\mu(f_1\otimes f_2)(z)=f_1(z) f_2(z)$,
\item the unit $\gone:\gC\to H$ by $\gone(\lambda)(z)= \lambda$,
\item the coproduct $\Delta:H\to H\widehat\otimes_\tau H$ by $\Delta f( z_1,z_2)=f(z_1+z_2)$,
\item the counit $\eps: H\to\gC$ by $\eps(f)=f(0)$,
\item the antipode $S:H\to H$ by $S f(z)=f(-z)$,
\end{itemize}
where $f_i\in H$, $z_i\in \gR^{m|n}$, $\lambda\in\gC$. These maps have been shown to be continuous for the Fr\'echet topologies introduced above, which means that $H$ is a graded commutative Fr\'echet-Hopf algebra. See \cite{Bonneau:2003vb} for structures of topological Hopf algebras associated to Lie groups.
\medskip

Let $\algA$ be a Fr\'echet algebra on which there exists an action $\rho$ of $\gR^{m|n}$ satisfying the axioms of subsection \ref{subsec-udf}. We denote by $\mu_0:\algA\widehat\otimes_\pi\algA\to\algA$ the (undeformed) product of $\algA$, which is continuous. We can associate to the universal deformation of Theorem \ref{thm-udf} a twist $F\in\caL(\algA^\infty\widehat{\otimes}_\pi\algA^\infty)$:
\begin{equation*}
F=\kappa\int_{\gR^{m|n}\times\gR^{m|n}}\dd z_1\dd z_2\ e^{-\frac{2i}{\theta}\tilde\omega(z_1,z_2)}\rho_{z_1}\otimes\rho_{z_2},
\end{equation*}
which is a continuous endomorphism and such that $\forall a,b\in\algA^\infty$, $\mu_0 F(a\otimes b)=a\star_\rho b$. The map $\mu_F:=\mu_0 F:\algA^\infty \widehat\otimes_\pi\algA^\infty\to\algA^\infty$ is the deformed product on the smooth vectors $\algA^\infty$.
\begin{proposition}[Comodule-algebras]
\label{prop-comod}
The set $(\algA^\infty,\mu_0)$ is a Fr\'echet comodule-algebra on $H$, with continuous coaction $\chi:\algA^\infty\to\algA^\infty\widehat{\otimes}_\tau H$ given by
\begin{equation*}
\forall a\in\algA^\infty,\,\forall z\in\gR^{m|n}\quad:\quad \chi(a)(z):=\rho_z(a).
\end{equation*}
The deformed algebra $(\algA^\infty,\star_\rho)$ corresponds then to twisting by $F$ this $H$-comodule-algebra. Moreover, $(\algA^\infty,\mu_F)$ is also a Fr\'echet comodule-algebra on $H$.
\end{proposition}
The fact that $(\algA^\infty,\mu_F)$ be still a comodule-algebra on $H$ after deformation is not automatic. The reason is that the supergroup $\gR^{m|n}$ is abelian, and it means that the external symmetries of the deformation are classical and coincide with the graded commutative Hopf algebra $H$.

\subsection{Operator superalgebras}
\label{subsec-oper}

For usual geometry, the structures of Hilbert space and C*-algebra are adapted to describe functional spaces like $L^2(M)$, and $C(M)$ or $L^\infty(M)$, for a compact manifold $M$. In the case of supergeometry, we have to introduce the notions of Hilbert superspace and C*-superalgebra \cite{Bieliavsky:2010su}. We will also see that the structure of C*-algebra is adapted to the universal deformation formula of subsection \ref{subsec-udf} while the one of C*-algebra is not.

\begin{definition}[Hilbert superspace]
Let $\ehH=\ehH_0\oplus\ehH_1$ be a complex $\gZ_2$-graded vector space endowed with a scalar product $(\fois,\fois)$ turning it into a Hilbert space satisfying $(\ehH_0,\ehH_1)=\algzero$.

$(\ehH,J)$ is called a Hilbert superspace of parity $n\in\gZ_2$ if $J\in\caL(\ehH)$ is a homogeneous unitary operator of degree $n$ such that
\begin{equation*}
\forall x\in\ehH\quad:\quad J^2(x)=(-1)^{(n+1)|x|} x
\end{equation*}
\end{definition}
The operator $J$ permits to define another (but superhermitian) scalar product:
\begin{equation*}
\forall x,y\in\ehH\quad:\quad \langle x,y\rangle=(J(x),y).
\end{equation*}
The direct sum of two Hilbert superspaces of parity $n$ is a Hilbert superspace of parity $n$ while the tensor product of two Hilbert superspaces of parity $n_1$ and $n_2$ is a Hilbert superspace of parity $n_1+n_2$. Note that there exists different notions of Hilbert superspaces (see \cite{Rudolph:2000} for example). This one is closely related to Krein spaces \cite{Bognar:1974}.

\begin{proposition}
The space $\caL(\ehH)$ of bounded operators on a Hilbert superspace $(\ehH,J)$ is a $\gZ_2$-graded C*-algebra (for the usual adjoint), endowed with the superinvolution\footnote{i.e. an antilinear map on $\caL(\ehH)$ of degree 0 satisfying $(T^\dag)^\dag=T$ and $(S T)^\dag=(-1)^{|S||T|}T^\dag S^\dag$.} defined by
\begin{equation}
\forall x,y\in\ehH\quad:\quad \langle T^\dag(x),y\rangle=(-1)^{|T||x|}\langle x,T(y)\rangle\label{eq-superinvol}
\end{equation}
and called the superadjoint.
\end{proposition}
Note that the superadjoint $T^\dag$ is different from the usual adjoint $T^\ast$ (with respect to the scalar product $(\fois,\fois)$): $T^\dag(x)=(-1)^{(n+1)(|T|+|x|)+|T||x|}JT^\ast J(x)$.

\begin{definition}[C*-superalgebra]
A C*-superalgebra $\algA$ is a complete $\gZ_2$-graded subalgebra of $\caL(\ehH)$, for a certain Hilbert superspace $\ehH$, closed for the superinvolution \eqref{eq-superinvol}.
\end{definition}
A C*-superalgebra is not necessarily closed under the (usual) involution of $\caL(\ehH)$, and therefore is not a C*-algebra in general. These structures are useful in supergeometry: if $\ehH=L^2(\gR^{m|n})$ and $J=\ast$ is the Hodge operator \eqref{eq-hodge}, then $(\ehH,J)$ is a Hilbert superspace of parity $n\text{ mod }2$, and $L^\infty(\gR^{m|n})$, endowed with the norm $\norm f\norm=\sum_I \norm f_I\norm_\infty$, is a C*-superalgebra.

\subsection{UDF for C*-superalgebras}
\label{subsec-udfcstar}

Let us see in this subsection how the UDF (Theorem \ref{thm-udf}) can be applied in the context of C*-superalgebras. We denote by $\algA$ a C*-superalgebra contained in $\caL(\ehH)$ for a Hilbert superspace $\ehH$, and we consider an action $\rho:\gR^{m|n}\times(\algA\otimes\superA)\to(\algA\otimes\superA)$ satisfying the axioms of the beginning of subsection \ref{subsec-udf} and the additional ones:
\begin{equation*}
\forall a\in\algA\quad:\quad
(\rho^a)^\dag=\rho^{(a^\dag)},\qquad |\rho^a|=|a|.
\end{equation*}
The expression \eqref{eq-omf} of the quantization map $\Omega$ can be directly extended to $\algA$-valued functions of $\caB_\algA(\gR^{m|n})$ and Proposition \ref{prop-qmap} insures its continuity:
\begin{equation*}
\Omega:\caB_\algA(\gR^{m|n})\to\caL(L^2(Q)\otimes\ehH),
\end{equation*}
since $\caL(L^2(Q))\otimes\algA\subset \caL(L^2(Q)\otimes\ehH)$. Using Equation \eqref{eq-rhobound}, we define $\Xi:\algA^\infty\to\caL(L^2(Q)\otimes\ehH)$ by $\Xi(a)=\Omega(\rho^a)$, for $a\in\algA^\infty$.

\begin{theorem}[UDF for C*-superalgebras]
\label{thm-udfcstar}
The map $\Xi:(\algA^\infty,\star_\rho)\to\caL(L^2(Q)\otimes\ehH)$ is an injective algebra-morphism of degree 0 compatible with superinvolutions. We then define the norm $\norm a\norm_\rho:=\norm\Xi(a)\norm$ on $\algA^\infty$ and $\algA_\rho$ its completion for this norm.

$(\algA_\rho,\star_\rho,{}^\dag,\norm\fois\norm_\rho)$ is a C*-superalgebra, the deformation of $\algA$.
\end{theorem}

\section{Applications}

\subsection{To supergeometry}
\label{subsec-supergeom}

We now apply Theorem \ref{thm-udfcstar} on supergeometrical examples. More precisely, we can deform compact trivial supermanifolds on which the supergroup $\gR^{m|n}$ acts. Note that a supermanifold $X$ is said to be trivial if there exists a global odd coordinates system $\{\xi^I\}$ on $X$, which means that Berezin integration, Hodge operation and the space $L^2(X)$ (see subsection \ref{subsec-setting}) are well-defined for $X$.
\medskip

Let $X = X_o \times \gR^{0|q}$ be a compact trivial supermanifold of dimension $p|q$ (i.e. $X_o$ is a supermanifold of dimension $p|0$ whose body $\gB X_o = \gB X$ is compact). We also consider a smooth action $\tau: \gR^{m|n}\times X\to X$ of the abelian supergroup $\gR^{m|n}$ on $X$. By denoting elements of $X$ as $u=(v,\eta)$, with $v\in X_o$ and $\eta\in\gR^{0|q}$, and elements  of $\gR^{m|n}$ as $z=(y,\xi)$, with $y\in \gR^{m|0}$ and $\xi \in \gR^{0|n}$, we can decompose $\tau_z(u) \in X_o \times \gR^{0|q}$ with respect to powers of the odd coordinates as follows: 
\begin{equation*}
\forall z=(y,\xi)\in M = \gR^{m|n}
\quad:\quad
\tau_z u=\left((\tau_y v)^0_{IJ}\xi^I\eta^J,(\tau_y v)^1_{IJ}\xi^I\eta^J\right),
\end{equation*}
with $(\tau_y v)^i_{IJ}$ smooth functions on $\gR^m\times \gB X$, such that $(\tau_y v)^0_{IJ}$ and $(\tau_y v)^1_{IJ}$ take values respectively in $\gB X$ and in $\gR^q$. Inspection of degrees gives that $(\tau_y v)^i_{IJ}=0$ for $i= 1+|I|+|J|$.

In the following, we also assume the additional technical conditions that $\forall (I,J)\neq (\emptyset,\emptyset)$ we have $(\tau_y v)^0_{IJ}=0$, and that every component $(\tau_y v)^{1,k}_{IJ}$ ($k\in\{1,\dots,q\}$) of $(\tau_y v)^1_{IJ}$ is uniformly bounded in $y$.

\begin{proposition}[Deformation of supermanifolds]
\label{prop-superman}
Let the space $\caC(X)$ be the completion of $C^\infty(X)$ for the norm $\norm f\norm=\sum_I\norm f_I\norm_\infty$. $\algA=\caC(X)$ is a C*-superalgebra. We define the natural action of $\gR^{m|n}$ on $\algA$ by: 
\begin{equation*}
\forall f\in\caC(X)\ ,\ \forall z\in \gR^{m|n}\ ,\ \forall u\in X
\quad:\quad \rho_z (f)(u)=f(\tau_{-z}u).
\end{equation*}
Then the axioms for $\rho$ of the beginning of subsections \ref{subsec-udf} and \ref{subsec-udfcstar} are satisfied, $C^\infty(X)\subset\algA^\infty$, and Theorem \ref{thm-udfcstar} applies. It yields a deformation (C*-superalgebra) $(\algA_\rho,\star_\rho,{}^\dag,\norm\fois\norm_\rho)$ of the supermanifold $X$.
\end{proposition}

\subsection{To QFT}
\label{subsec-qft}

We now apply the above construction to provide an interpretation of a certain renormalizable quantum field theory (QFT) on the Moyal space (see \cite{Wulkenhaar:2006si,deGoursac:2009gh} for a review), which is the deformation quantization of $\gR^m$. This has been done in \cite{deGoursac:2010zb,Bieliavsky:2010su}.
\medskip

The usual $\phi^4$ real scalar field theory on the Euclidean Moyal space $\gR^m_\theta$, with parameters $\nu$ and $\lambda$,
\begin{equation*}
S(\phi)=\int \dd^mx\Big(\frac 12(\partial_\mu\phi)^2 +\frac{\nu^2}{2}\phi^2 +\lambda\phi\star\phi\star\phi\star\phi\Big)
\end{equation*}
suffers from the Ultraviolet-Infrared mixing \cite{Minwalla:1999px}, a new type of divergence generic to noncommutative spaces \cite{Gayral:2005af} which spoils its renormalizability, a very important property of QFT. Note also that the action can be reformulated by using $\partial_\mu\phi=[-\frac i2\wx_\mu,\phi]_\star$ where $\wx=\frac{2}{\theta}\omega_0(x,\fois)$. This problem of UV-IR mixing has been solved by the addition of a harmonic term in the action:
\begin{equation}
S(\phi)=\int \dd^mx\Big(\frac 12[-\frac i2\wx_\mu,\phi]_\star^2 +\frac{\Omega^2}{2}\wx^2\phi^2+\frac{\nu^2}{2}\phi^2 +\lambda\phi\star\phi\star\phi\star\phi\Big),\label{eq-actharm}
\end{equation}
where $\Omega$ is a non-zero parameter. The resulting theory is renormalizable to all orders \cite{Grosse:2004yu} for $m=2,4$, and possesses new interesting properties \cite{Disertori:2006nq,deGoursac:2007uv,Tanasa:2007xa,deGoursac:2009fm,Hounkonnou:2009qt}. Note that there exists now another renormalizable real scalar theory on the Moyal space \cite{Gurau:2008vd,Tanasa:2010fk}. A gauge theory associated to \eqref{eq-actharm} has been constructed in \cite{deGoursac:2007gq,Grosse:2007dm}. It exhibits similar features \cite{deGoursac:2008rb,Cagnache:2008tz} and is a candidate to renormalizability. See \cite{Blaschke:2007vc,Blaschke:2009aw} for the study of the BRST symmetry.
\medskip

We want to give an interpretation to the harmonic term in \eqref{eq-actharm}. Let us consider $\caB(\gR^{m|1})$ endowed with the star-product \eqref{eq-moy}. We denote by $[-,-]_\star$ the graded star-bracket of superfunctions\footnote{do not confuse with the Moyal bracket, it is the same notation.} and by $\Tr$ the trace $\Tr(f)=\int \dd^m x\ f(x,0)$ where $f(x,\xi)$ is a superfunction on $\gR^{m|1}$. Let $\phi$ be as above a function on $\gR^m$ and $b\in \gR$. Then, $(x,\xi)\mapsto(1+b\xi)\phi(x)$ is a superfunction on $\gR^{m|1}$ and the standard $\phi^4$ action for this field:
\begin{equation*}
\Tr\Big(\frac12|[-\frac i2(1+b\xi)\wx_\mu,(1+b\xi)\phi]|^2+\frac{\nu^2}{2}|(1+b\xi)\phi|^2 +\Lambda |((1+b\xi)\phi)\star((1+b\xi)\phi)|^2\Big)
\end{equation*}
coincides exactly with \eqref{eq-actharm} for $\Omega=\frac{\alpha\theta b^2}{(1+\alpha)^2}$ and $\lambda=\Lambda(1+\Omega^2)$. It shows that the noncommutative supergeometry of the deformation of $\gR^{m|1}$ is likely to be related to the renormalizability of the noncommutative scalar field theory \eqref{eq-actharm}. This interpretation could be useful to obtain other renormalizable noncommutative quantum field theories.

\vskip 1 true cm

{\bf Acknowledgements}: The author thanks the organizers of the EU-NCG 4th Annual Meeting (Bucharest, April 2011).

\bibliographystyle{utcaps}
\bibliography{biblio-these,biblio-perso,biblio-recents}

\end{document}